\documentclass[a4paper]{amsart}

\usepackage[english]{babel}
\usepackage[utf8]{inputenc}
\usepackage[T1]{fontenc}
\usepackage{lmodern}

\usepackage{natbib}

\usepackage{graphicx}

\usepackage{amssymb} 
\usepackage{version} 
\newcommand{\defegal}{\triangleq}
\newcommand{\espacea}[1]{\mathbb{#1}}
\newcommand{\va}[1]{\boldsymbol{\uppercase{#1}}}
\newcommand{\esper}[2][]{\mathbb{E}_{#1}\left(#2\right)}
\newcommand{\espcond}[3][]{\mathbb{E}_{#1}\left(\left. #2 \right\vert #3\right)}
\newcommand{\omeg}{\Omega}
\newcommand{\trib}{\mathcal{A}}
\newcommand{\prbt}{\mathbb{P}}
\newcommand{\bbR}{\mathbb{R}}
\newcommand{\findi}[1]{\mathbf{1}_{#1}}

\theoremstyle{remark}
\newtheorem{remark}{Remark}
\newtheorem{hyp}{Assumption}

\usepackage[hypertexnames=false,colorlinks=true,linkcolor=black,citecolor=black]{hyperref}

\begin{document}

\title{A comparison of sample-based Stochastic Optimal Control methods}

\author[P. Girardeau]{Pierre Girardeau}
\address{P. Girardeau, EDF R\&D, 1, avenue du Général de Gaulle, F-92141 Clamart Cedex, France, also with Universit\'e Paris-Est, CERMICS and ENSTA.}
\email{pierre.girardeau@cermics.enpc.fr}


\date{\today}


\thanks{The author would like to thank his PhD supervisors Pierre Carpentier and Guy Cohen for their helpful remarks on previous versions of this paper.}

\begin{abstract}
	In this paper, we compare the performance of two scenario-based numerical methods to solve stochastic optimal control problems: scenario trees and particles. The problem consists in finding strategies to control a dynamical system perturbed by exogenous noises so as to minimize some expected cost along a discrete and finite time horizon. We introduce the Mean Squared Error (MSE) which is the expected $L^2$-distance between the strategy given by the algorithm and the optimal strategy, as a performance indicator for the two models. We study the behaviour of the MSE with respect to the number of scenarios used for discretization. The first model, widely studied in the Stochastic Programming community, consists in approximating the noise diffusion using a scenario tree representation. On a numerical example, we observe that the number of scenarios needed to obtain a given precision grows exponentially with the time horizon. In that sense, our conclusion on scenario trees is equivalent to the one in the work by~\citet{Shapiro06} and has been widely noticed by practitioners. However, in the second part, we show using the same example that, by mixing Stochastic Programming and Dynamic Programming ideas, the particle method described by~\citet{CarpentierCohenDallagi09} copes with this numerical difficulty: the number of scenarios needed to obtain a given precision now does not depend on the time horizon. Unfortunately, we also observe that serious obstacles still arise from the system state space dimension.
\end{abstract}

\maketitle

\section*{Introduction}

Consider a controlled dynamical system, affected by some exogenous noises, say uncertain parameters for instance. Stochastic optimal control consists in driving this system, having at each time step partial or total observations of those noises, so as to minimize some expected cost integrated through time and while satisfying a number of constraints. Many applications can be handled through such a model. For instance, think of a power producer that has to plan the use of a set of power plants (the control would be the production of the plants at each time step) in order to supply an uncertain power demand (the noises) while minimizing a production cost over a certain time period. We here consider discrete and finite time horizon. We are hence looking for feedback functions which, at each time step and for each possible observation of the system, provide a decision to be taken. We consider the laws of the random variables involved to be continuous; consequently this is an infinite-dimensional optimization problem. In most situations, no analytic solution can be found and one has to consider approximations of the original problem to be able to solve it numerically.

In the manner of the Monte Carlo approach for computing the expectation of a random variable, we seek a tractable approximation of the probabilistic structure of the original problem based on samples of some noise scenarios. A scenario consists in a sample of the noises along the whole time horizon. Note that such a resolution method is stochastic in itself. Indeed, the solution given by such a method is a random variable in the sense that it depends on the scenarios sampled and used in the resolution. Therefore, when studying the performance of an algorithm based on such an approach, a variance term naturally appears, representing the sensitivity of the solution regarding these scenarios. This variance term is added to the squared bias term which arises from the necessity to approximate the solution in feedback by a function in a predefined space. Note that it is also present in the performance evaluation of deterministic techniques, like Dynamic Programming~\citep{BertsekasDP} for instance.

For stochastic optimal control problems, it is common to represent the diffusion of ``likely futures'' using a scenario tree structure, leading to so-called multi-stage stochastic programs. This kind of representation goes back to~\citet{Dantzig55} and has been widely studied within the Stochastic Programming community~\citep[see][for a broad overview of this approach]{PrekopaStoProg,StochasticProgramming03}. This methodology consists in discretizing the probabilistic structure of the problem, then rewriting the constraints and objective function of the original problem on this discrete structure and finally solving it using a well-suited mathematical programming technique. In other words, the problem is solved for a particular sample of the random variables (a scenario tree) and the solution is hence a random variable in itself. We have to keep this fact in mind while evaluating the error. The main interest of such a methodology is that it leads back to a deterministic problem, on which we can apply classical tools of numerical optimization~\citep{BonnansGilbertLemarechalSagastizabal06}. For instance, one may then try to solve large-scale problems using decomposition techniques~\citep[see][Ch. 3]{CarpentierCohenCulioli95,StochasticDecomposition96,StochasticProgramming03}.

Despite the benefits of this methodology, it has been already observed that, in order to obtain a given precision, the number of scenarios needed to build the scenario tree has to grow exponentially with the time horizon of the problem. One of the aims of this paper is to highlight and to quantify this practical observation on an example. Thus, we start in Section~\ref{sec:MathForm} by introducing the Mean Squared Error, that is the distance between the strategy obtained using scenario trees and the (supposedly unique) optimal strategy. This is the indicator we use to estimate the performance of the method. Then, in Section~\ref{sec:ScenarioTrees}, we study on a numerical example the relation linking the Mean Squared Error to the number of scenarios used for the discretization. \citet{Shapiro06} obtains similar conclusions by working on the performance regarding the objective function and using large deviations tools.

We show that this negative observation is to be attributed to the scenario tree approach rather than to be considered as a feature of the original problem. Indeed, when using the particle method~\citep{TheseDallagi,CarpentierCohenDallagi09}, in Section~\ref{sec:Particles} on the same example, we show that the number of scenarios needed to obtain a given precision does not grow when the time horizon gets longer. This methodology, which mixes ideas from Stochastic Programming and Dynamic Programming, is a variational technique which consists in writing the optimality conditions of the optimization problem first, and then solving them using sampling techniques. This is a gradient-like method that builds an adaptive mesh over the state space as gradient iterations progress towards the solution. This mesh aims at discretizing the space in the regions mostly visited by the state vector at the optimum. According to the results obtained on the example we present here, the algorithm seems to be well-suited for multi-stage problems. However, with this approach, it seems that the difficulties arising from the dimension of the state space remain; this point is discussed in Section~\ref{sec:Dimension}.

\section{Mathematical formulation} \label{sec:MathForm}

Consider a controlled dynamic system, affected by random variables called noises. We aim to find strategies that minimize some expected cost over a finite time horizon, while satisfying a number of constraints. We suppose that the problem has a unique solution and propose to compare the approximate strategies using the Mean Squared Error (MSE).

\subsection{The problem}

We denote by~$T$ the finite time horizon of the problem and consider discrete time~$\{0, \dots, T\}$. Let $(\omeg, \trib, \prbt)$ be a probability space. Three kinds of random variables on this space are involved in the problem\footnote{Throughout this paper, random variables will be denoted by bold letters (e.g. $\va{w} \in L^2(\omeg,\trib,\prbt;\bbR^p$).}:
\begin{itemize}
	\item the state\footnote{In the Stochastic Optimal Control framework, the terminology ``state'' has a rather precise sense: this is the minimal information that completely sums up the past of the system in order to compute its future behaviour and cost function knowing all future inputs~\citep[see][Def.~5.4.1]{Powell07}. In that sense, the variable~$\va{x}$ can only be considered as the state variable when assumption~\ref{hyp:Indpdt} is made later on in~\S\ref{ssec:MSE}. Until then, our use of this terminology is somewhat abusive.}~$\va{x} = (\va{x}_t)_{t=0, \dots, T}$, the values of which lie in a finite-dimensional vector space~$\espacea{X}$;
	\item the control~$\va{u} = (\va{u}_t)_{t=0, \dots, T-1}$, the values of which lie in a finite-dimen\-sional vector space~$\espacea{U}$;
	\item the noise~$\va{w} = (\va{w}_t)_{t=0, \dots, T}$, the values of which lie in a finite-dimensional vector space~$\espacea{W}$, equipped with the~$\sigma$-field~$\mathcal{W}$ and the probability~$\prbt_{\va{w}}$, which is the transport of the probability~$\prbt$ by~$\va{w}$.
\end{itemize}
Since some of the numerical methods we use in the following are gradient-based, it is natural to assume that the state and control lie in a Hilbert space. Hence we assume all three random variables are in~$L^2(\omeg,\trib,\prbt)$. The noise variable~$\va{w}$ is given, i.e. its probability law is known, while the state and control variables are optimization variables. We here suppose that the information structure has perfect memory: at each time step~$t$, the noise~$\va{w}_t$ is observed and kept in memory, in such a way that the decision at time~$t$ is based on the knowledge of~$(\va{w}_0, \dots, \va{w}_t)$. Hence we introduce the $\sigma$-field~$\mathcal{F}_t = \sigma\{\va{w}_0, \dots, \va{w}_t\}$ that represents the information available at time~$t$, the associated filtration~$(\mathcal{F}_t)_{t=0, \dots, T}$, and require the control at time~$t$ to be measurable with respect to~$\mathcal{F}_t$. This constraint will be written:~$\va{u}_t \preceq \mathcal{F}_t$.

At time~$0$, we assume that the value of the state is the realization of some random variable~$\va{w}_0$. Then, at each time step $t$, based on the available information, a decision~$\va{u}_t$ is taken, a cost~$C_t(\va{x}_t, \va{u}_t)$ is incurred and the state evolves according to~$\va{x}_{t+1} = f_t(\va{x}_t,\va{u}_t, \va{w}_{t+1})$. At the end of the time period, a final cost~$V(\va{x}_T)$ is incurred.

The aim is to find a pair~$(\va{x}, \va{u})$ that minimizes the expected sum of the costs over the time horizon while satisfying the dynamics and measurability constraints. Mathematically speaking, the problem we consider may be formulated as follows:
\begin{subequations} \label{eqn:P}
\begin{align}
	\min_{\va{x}, \va{u}} \quad& \esper{\sum_{t=0}^{T-1} C_t\left(\va{x}_t, \va{u}_t\right) + V\left(\va{x}_T\right)}, \label{eqn:P-1} \\
	\text{s.t.} \quad& \va{x}_{t+1} = f_t\left(\va{x}_t, \va{u}_t, \va{w}_{t+1}\right), \qquad \forall t=0, \dots, T-1, \label{eqn:P-2} \\
	& \va{x}_0 = \va{w}_0, \label{eqn:P-3} \\
	& \va{u}_t \preceq \mathcal{F}_t, \qquad \forall t=0, \dots, T. \label{eqn:P-4}
\end{align}
\end{subequations}
Equations~\eqref{eqn:P-2} and \eqref{eqn:P-3} describe the dynamics of the state variable and are~$\prbt$-almost sure relations. Equation~\eqref{eqn:P-4} states the information structure of the problem: it is the measurability constraint that forces the decision to be a function of all the past noises, and to be independent of future noise values. Hence we refer to this relation as the non-anticipativity constraint.

\subsection{Mean Squared Error (MSE) for strategies} \label{ssec:MSE}

Now suppose that Problem~\eqref{eqn:P} has a unique solution denoted by~$(\va{x}^*, \va{u}^*)$. Suppose that we also have a numerical method, based on a scenario discretization, that gives us an approximate solution of~\eqref{eqn:P} that we denote by~$(\va{x}^\sharp, \va{u}^\sharp)$. Since this approximate solution given by the numerical method depends on, say,~$N$ samples, the pair~$(\va{x}^\sharp, \va{u}^\sharp)$ lies in the probability space associated with these samples, that is~$(\espacea{W}, \mathcal{W}, \prbt_{\va{w}})^{\otimes N}$. For instance, scenario tree-based methods compute values of the optimal decision for a finite number of state values that depend on the samples used for building the tree structure.

We make the following assumption:
\begin{hyp} \label{hyp:Indpdt}
	Random variables~$\va{w}_0, \dots, \va{w}_T$ are independent.
\end{hyp}
Under this assumption, the problem lies in the framework of Dynamic Programming~(DP). In particular, we know that optimal controls can be expressed as feedback functions depending only on the state variable~$x$. In other words, there exists some sequence of functions~$\gamma^*$ such that~$\va{u}_t^*=\gamma_t^*(\va{x}_t^*)$.

In order to be able to compare the optimal and approximate solutions, we need to produce a solution of the same nature as the optimal solution, that is a feedback function such as~$\gamma^*$ out of~$(\va{x}^\sharp, \va{u}^\sharp)$. We hence suppose that we have such an interpolation/regression operator that gives us a feedback function~$\va{\Gamma}^\sharp = (\va{\Gamma}_t^\sharp)_{t=0, \dots, T}$ out of the approximate solution~$(\va{x}^\sharp, \va{u}^\sharp)$. More details on the interpolation/regression operators we use are given in \S\ref{ssec:TreesCaseStudy}. Note that the approximate strategy~$\va{\Gamma}^\sharp = (\va{\Gamma}_t^\sharp)_{t=0, \dots, T-1}$ is a random variable that lies in the samples probability space introduced earlier, that is~$(\espacea{W}, \mathcal{W}, \prbt_{\va{w}})^{\otimes N}$.

We are interested in evaluating the efficiency of our numerical method for solving Problem~\eqref{eqn:P} by computing the distance between the approximate strategy $\va{\Gamma}^\sharp$ and the optimal one $\gamma^*$. Since those functions take values over the state space, we need to define some measure on this space. In our context, the measure that seems the most natural is the one corresponding to the density of the optimal state~$\va{x}^*$. Indeed, using this measure, we allocate more weight to regions which are often ``visited'' by the optimal state and we do not take into account decisions in regions that are never ``visited'' by the optimal state, since those decisions are unlikely to be used in practice. We denote by~$\mu_t^*$ the density of the optimal state at time~$t$.

We can now introduce the MSE associated with strategy~$\va{\Gamma}^\sharp$, which is our performance indicator in this study:
\begin{equation*}
	\text{MSE} \defegal \esper{\sum_{t=0}^{T-1} \int_{\espacea{X}} \left\Vert \gamma_t^*\left(x\right) - \va{\Gamma}_t^\sharp\left(x\right) \right\Vert^2 \mu_t^*\left(x\right) \text{d}x}.
\end{equation*}
Note that the expectation in the MSE is taken over all possible samples; in other words, it lies on the samples probability space~$(\espacea{W}, \mathcal{W}, \prbt_{\va{w}})^{\otimes N}$. The MSE is thus the expected distance between the optimal strategy and the approximate strategy. Recall that this expectation is taken with respect to the drawings used by the algorithm to produce strategy~$\va{\Gamma}^\sharp$. It is common to decompose the MSE using the so-called variance and squared bias. For this purpose, we introduce the expectation~$\gamma^\sharp$ of the approximate feedback function~$\va{\Gamma}^\sharp$. Then
\begin{align}
	\text{MSE} &= \underbrace{\sum_{t=0}^{T-1} \int_{\espacea{X}} \left\Vert \gamma_t^*\left(x\right) - \gamma_t^\sharp\left(x\right) \right\Vert^2 \mu_t^*\left(x\right) \text{d}x}_{\text{Squared bias}} \nonumber \nonumber \\
	&\qquad + \underbrace{\esper{\sum_{t=0}^{T-1} \int_{\espacea{X}} \left\Vert \gamma_t^\sharp\left(x\right) - \va{\Gamma}_t^\sharp\left(x\right) \right\Vert^2 \mu_t^*\left(x\right) \text{d}x}}_{\text{Variance}}. \label{eqn:MSE}
\end{align}
The squared bias is the square of the distance between the optimal strategy and the expected approximate strategy, while the variance is the expected square of the distance between the approximate strategy and the expected approximate strategy. Both are real values.

\begin{remark}
	Another relevant choice as a performance indicator would have been the following:
	\begin{enumerate}
		\item Having the approximate strategy~$\va{\Gamma}^\sharp$, build the``true'' state and control obtained while simulating this strategy, that is the pair~$(\va{x}^\dag, \va{u}^\dag)$ that satisfies the dynamics~\eqref{eqn:P-2} and~\eqref{eqn:P-3}:
		\begin{align*}
			\va{x}_0^\dag &\defegal \va{w}_0, \\
			\va{x}_{t+1}^\dag &\defegal f_t\left(\va{x}_t^\dag, \va{u}_t^\dag, \va{w}_{t+1}\right) \text{, with } \va{u}_t^\dag \defegal \va{\Gamma}_t^\sharp(\va{x}_t^\dag), \qquad \forall t=0, \dots, T-1.
		\end{align*}
		\item Since the optimal decision~$\va{u}^*$ is a random variable on~$(\omeg, \trib, \prbt)$ and the decision~$\va{u}^\dag$ is a random variable that lies on the tensor product of~$(\omeg, \trib, \prbt)$ and the samples probability space~$(\espacea{W}, \mathcal{W}, \prbt_{\va{w}})^{\otimes N}$, compute the expectation:
		\begin{equation*}
			\esper{\left\Vert \va{u}^\dag - \va{u}^* \right\Vert_{\espacea{U}}^2},
		\end{equation*}
		where the expectation lies on the tensor product of the probability spa\-ces~$(\omeg, \trib, \prbt)$ and~$(\espacea{W}, \mathcal{W}, \prbt_{\va{w}})^{\otimes N}$. Note that this quantity may also be written~$\mathbb{E}(\Vert \va{u}^\dag - \va{u}^* \Vert_{L^2\left(\omeg, \trib, \prbt\right)}^2)$, where the expectation now lies on the samples probability space only.
	\end{enumerate}
	This quantity is another relevant performance indicator since it is positive and it equals zero if and only if~$\va{u}^\dag$ equals~$\va{u}^*$ almost surely. Moreover, we note that the choice of~$\va{u}^\dag$ for evaluating the performance of the numerical method is the one made by~\citet{Shapiro06}, even though the comparison is performed using the objective function rather than the strategies.
\end{remark}

\subsection{Computing the MSE} \label{ssec:CompMSE}

For computing the squared bias and the variance as defined in Equation~\eqref{eqn:MSE}, we must face two main issues.
\begin{enumerate}
	\item At each time step~$t$, we must compute integrals over the whole state space with respect to the density~$\mu_t^*$. In the following example, we are able to explicitly compute the optimal control strategy. We can then numerically integrate the Fokker-Planck equation to derive the density of the optimal state on a sufficiently dense grid and perform the integration by quadrature techniques. However, this technique suffers for the same \textit{curse of dimensionality} as DP. For this reason we choose to perform the integration using a Quasi-Monte Carlo technique. This sampling method aims at distributing the samples associated with density~$\mu_t^*$ equally over the state space using so-called low-discrepancy quasi-random sequences. This provides a very efficient convergence speed for the numerical computation of expectations. There exist numerous methods of Quasi-Monte Carlo sampling~\citep[see][for a survey of these methods]{Niederreiter92}.
	\item We must be able to compute the expectation in the variance term. This expectation is taken over the samples used during the algorithm to produce the approximate strategy, e.g. the scenarios in the case of a scenario tree technique. We approximate this expectation by Monte Carlo, by performing the experiment a great number of times independently one from another. In our case,~$10^4$ experiments were enough to precisely evaluate the variance.
\end{enumerate}

We are now able to estimate the efficiency of numerical methods using the MSE. We point out the relation between the error and the parameters of the methods, namely the number of scenarios used to discretize randomness.

\section{Scenario tree-based methods} \label{sec:ScenarioTrees}

In the context of discrete time and finite horizon stochastic optimal control problems, a popular and widely used resolution technique consists in approximating the information structure of Problem~\eqref{eqn:P} using scenario trees. This methodology has been widely studied by the Stochastic Programming community. We do not detail this method and refer to~\citet{StochasticProgramming03} for a more in depth presentation of Stochastic Programming.

\subsection{Brief presentation}

Because of the measurability constraints in Problem~\eqref{eqn:P}, and since there is no reason for the noise probability densities to have finite support, decision variables~$\va{x}$ and~$\va{u}$ are in general infinite-dimensional. Hence, except for some very particular cases where an explicit solution can be found, solving the underlying optimization problem directly is intractable. Scenario trees provide a way to approximate filtrations using noise samples to produce finite support approximations.

Let~$\mathcal{N}$ be the set of all nodes in the tree,~$\mathcal{R}$ be the set of the root nodes and~$\mathcal{L}$ the set of the leaves. Building a scenario tree consists in defining the following functions:
\begin{itemize}
	\item the time function~$\theta: \mathcal{N} \rightarrow \{0, \dots, T\}$ which, with every node, associates the corresponding time step and the corresponding multi-application~$\theta^{-1}$ which, to every time step~$t \in \{0, \dots, T\}$, associates the set of all nodes at time~$t$;
	\item the function~$\nu: \mathcal{N} \backslash \mathcal{R} \rightarrow \mathcal{N} \backslash \mathcal{L}$ which, with each non-root node, associates the preceding node in the tree;
	\item the weight function~$\pi: \mathcal{N} \rightarrow [0,1]$ which, with every node in the tree, associates the probability to go through that node (the sum of all~$\pi(i)$, with $i \in \theta^{-1}(t)$ is equal to $1$, for each $t$);
	\item the multi-application~$F=\nu^{-1}$ which, with every node~$i$ in the tree, associates the set of its successors (we impose the convention that~$F(i) = \emptyset, \forall i \in \mathcal{L}$);
	\item the function~$F^+$ which, with each node~$i$, associates the whole sub-tree starting from~$i$: $F^+(i)=F(i) \cup F^2(i) \cup \dots \cup F^{T-\theta(i)}(i)$.
\end{itemize}
We consider a regularly branching tree, i.e. with a constant branching factor denoted by~$n_b$, so that~$\text{card}(F(i))=n_b$, for every node $i \in \mathcal{N} \backslash \mathcal{L}$. At time~$t=0$, we draw a~$n_b$-sample of the random variable~$\va{w}_0$, that we denote by~$(\va{w}^{\prime i})_{i \in \theta^{-1}(0)}$. For each root node~$i$, we draw at time~$t=1$ a~$n_b$-sample of the random variable~$\va{w}_1$\footnote{This tree building procedure, often called conditional sampling, usually requires to draw samples of~$\va{w}_1$ knowing the value of the preceding noises, that is~$\va{w}_0$. Because of Assumption~\ref{hyp:Indpdt}, this is not needed here.}, that we denote by~$(\va{w}^{\prime j})_{j \in F(i)}$. Samples of the random variable~$\va{w}_1$ are drawn independently root node per root node. This procedure continues until final time~$t=T$, so that we have~$n_b^{T+1}$ leaves to the tree. Hence one can note that this can reveal a heavy procedure when~$T$ becomes large, since we need to draw~$n_b^{T+1}$ samples of~$\va{w}_T$.

With this material we can now define a discretized version of Problem~\eqref{eqn:P}:
\begin{subequations} \label{eqn:PDiscret}
\begin{align}
	\min_{\va{X}^\prime, \va{u}^\prime} \quad& \sum_{i \in \mathcal{N} \backslash \mathcal{L}} \pi(i) \cdot C_{\theta(i)}\left(\va{x}^{\prime i}, \va{u}^{\prime i}\right) + \sum_{i \in \mathcal{L}} V\left(\va{x}^{\prime i}\right), \label{eqn:PDiscret-1} \\
	\text{s.c.} \quad& \va{x}^{\prime i} = f_{\theta(\nu(i))}\left(\va{x}^{\prime \nu(i)}, \va{u}^{\prime \nu(i)}, \va{w}^{\prime i}\right), \qquad \forall i \in \mathcal{N} \backslash \mathcal{R}, \label{eqn:PDiscret-2} \\
	& \va{x}^{\prime i} = \va{w}^{\prime i}, \qquad \forall i \in \mathcal{R}. \label{eqn:PDiscret-3}
\end{align}
\end{subequations}
The expectation in the objective function~\eqref{eqn:P-1} has been replaced in~\eqref{eqn:PDiscret-1} by a discrete weighted sum of costs on the tree nodes. Equations~\eqref{eqn:PDiscret-2} and~\eqref{eqn:PDiscret-3} correspond to Equations~\eqref{eqn:P-2} and~\eqref{eqn:P-3}. The only constraint that seems to be missing is Equation~\eqref{eqn:P-4}. Actually, it is now reflected in the structure of the tree. Indeed, starting from any node in the tree, there is only one possibility for going back to the root (we know exactly the values of all past noises), but there is a whole sub-tree that goes to the set of leaves, corresponding to the final time step (we only know the laws of future noises). In other words, the non-anticipativity constraint is coded in the structure of the tree.

In most cases, the scenario tree is not drawn directly. One draws a number of scenarios independently from one another, and then one builds a tree structure. Building such a tree in an efficient way (regarding the implied discretization error) is a challenging task. There exists a large literature on the way one can build scenario trees in a reasonable way~\citep{Pflug,HeitschRomisch03}, as well as on stability of the underlying optimization problem regarding this discretization~\citep{HeitschRomischStrugarek05}. We will not focus on this point here. As it has already been mentioned, we suppose that the tree has been built using a branching factor~$n_b$, which does not depend on time, leading to~$n_b^{T+1}$ leaves. Recall that this procedure requires the drawing of~$N \defegal n_b^{T+1}$ different scenarios.

\subsection{Case study} \label{ssec:TreesCaseStudy}

We apply this methodology to a simple instance where Problems~\eqref{eqn:P} and~\eqref{eqn:PDiscret} can be solved analytically. Then, we apply the methodology described in \S\ref{ssec:CompMSE} in order to compute the MSE and to point out its relation with the number of scenarios~$N$ used to discretize the original problem. Let~$\varepsilon$ be some positive real number, the problem we consider is the following:
\begin{subequations} \label{eqn:example}
\begin{align}
	\min_{\va{x}, \va{u}}\quad& \esper{\varepsilon \sum_{t=0}^{T-1} \va{u}_t^2 + \va{x}_T^2}, \\
	\text{s.c.}\quad& \va{x}_{t+1} = \va{x}_t + \va{u}_t  + \va{w}_{t+1}, \quad \forall t=0, \dots, T-1, \\
	& \va{x}_0 = \va{w}_0, \label{eqn:example-3} \\
	& \va{u}_t \preceq \mathcal{F}_t.
\end{align}
\end{subequations}
The state and control variables lie in~$L^2(\omeg,\trib,\prbt;\bbR)$. Noise variables~$\va{w}_0$, $\dots$, $\va{w}_T$ are i.i.d. random variables with uniform law on~$[-1,1]$. One easily shows that the optimal control for Problem~\eqref{eqn:example} reads~$\va{u}_t^*=\gamma_t^*(\va{x}_t^*)$ with:
\begin{equation} \label{eqn:sol}
	\gamma_t^*\left(x\right) = -\frac{x}{T - t + \varepsilon}, \qquad \forall x \in \bbR.
\end{equation}

We now solve the problem on the tree. On each node~$i$, one has a noise sample~$\va{w}^{\prime i}$ and solving the problem on the tree leads to values~$\va{x}^{\prime i}$ and~$\va{u}^{\prime i}$ for the state and control, respectively. One easily shows that:
\begin{equation*}
	\va{u}^{\prime i} = -\frac{\va{x}^{\prime i} + \frac{1}{n_b^{T-\theta(i)}} \sum_{j \in F^+(i)} \va{w}^{\prime j}}{T-\theta(i)+\varepsilon}, \qquad \forall i \in \mathcal{N} \backslash \mathcal{L}.
\end{equation*}
We now need to build a feedback function from these state and control values. This can be performed using a interpolation/regression operator. We here use the simplest one, which is the nearest neighbour interpolation operator~\citep[see][\S10.4, for a rigorous study of this operator and its possible implementations]{GershoGray92}. At each time step~$t$, we build the Voronoi diagram of the set of state values~$\{\va{x}^{\prime i}\}_{i \in \theta^{-1}(t)}$. We denote by~$\va{c}^{\prime i}$ the cell for which~$\va{x}^{\prime i}$ is the center. This cell is a random variable in itself, because it depends on the random variables drawn to build the scenario tree. On this cell, the feedback function will be constant, equal to~$\va{u}^{\prime i}$. Thus we obtain the following strategy at time t:
\begin{equation} \label{eqn:Strategy}
	\va{\Gamma}_t^\sharp\left(x\right) = \sum_{i \in \theta^{-1}(t)} \va{u}^{\prime i} \findi{\va{c}^{\prime i}}\left(x\right).
\end{equation}

\begin{remark}
The fact that the Voronoi cells~$\va{c}^{\prime i}$ are random variables, depending on the scenarios, makes the theoretical study of the error associated with strategy~$\va{\Gamma}^\sharp$ intricate. Indeed, in order to compute the expected approximate strategy~$\gamma^\sharp$, as well as to evaluate the variance term in Equation~\eqref{eqn:MSE}, we have to compute expectations over all possible scenario drawings, leading through Equation~\eqref{eqn:Strategy} to expectations over all possible Voronoi cells. This is generally a challenging task. This is the reason why we here appeal to numerical experiments to study the behaviour of the error with respect to the number of scenarios.
\end{remark}

\begin{figure}[htbp]
\begin{center}
\begin{tabular}{cc}
\includegraphics[width=0.45\textwidth]{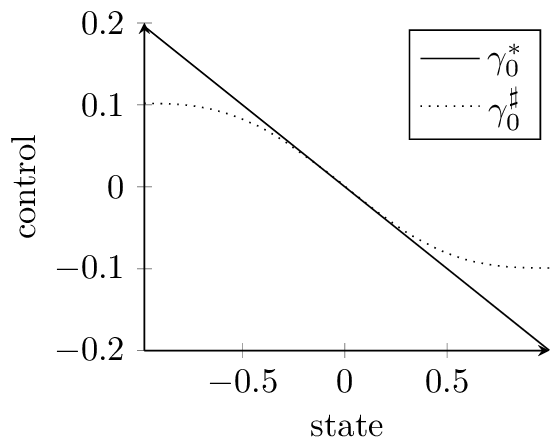}
&
\includegraphics[width=0.45\textwidth]{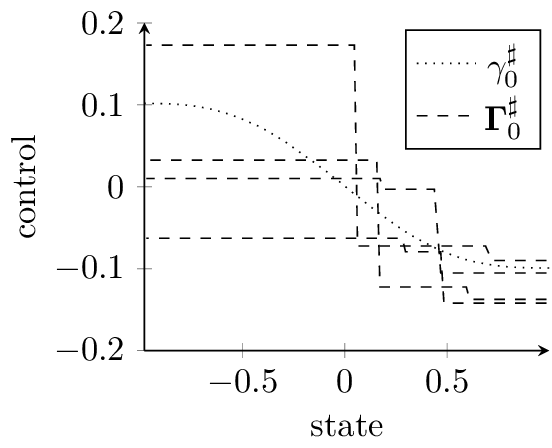}
\end{tabular}
\end{center}
\caption{\label{fig:StratExampleTrees}Strategies for scenario trees}
\end{figure}

The results for this approach are presented in Figure~\ref{fig:StratExampleTrees} for~$T=4$ and a branching factor of~$n_b=3$. The approximate strategy~$\va{\Gamma}_0^\sharp$ for the first time step is hence a piecewise constant function with~$3$ pieces. On the right-hand side of Figure~\ref{fig:StratExampleTrees}, we draw~$4$ samples of this strategy (dashed curves) to emphasize that the scenario-tree method is stochastic: each sample~$\va{\Gamma}_0^\sharp$ is derived from a different scenario tree. We also draw the averaged strategy~$\gamma_0^\sharp$ over the~$10^4$ scenarios (dotted curve). One can remark that even if sample strategies are not continuous functions, the average strategy seems to be continuous, even smooth. On the left-hand side of figure~\ref{fig:StratExampleTrees}, we draw the exact strategy (solid curve) obtained using equation~\eqref{eqn:sol} and the same average approximate strategy as on the right-hand side (dotted curve).

We can observe both the variance and the bias terms of the MSE on Figure~\ref{fig:StratExampleTrees}: the variance term is the average $L^2$-distance between the dashed curves~(of course if there were~$10^4$ of them) and the dotted curve on the right hand-side, whereas the bias term is the $L^2$-distance between the solid curve and the dotted curve on the left hand-side. From Equation~\eqref{eqn:example-3}, the density of the optimal state, under which we compute the $L^2$-distances cited above, equals the density of~$\va{w}_0$, which follows a uniform law on~$[-1,1]$.

\begin{figure}[htbp]
\begin{center}
\begin{tabular}{cc}
\includegraphics[width=0.45\textwidth]{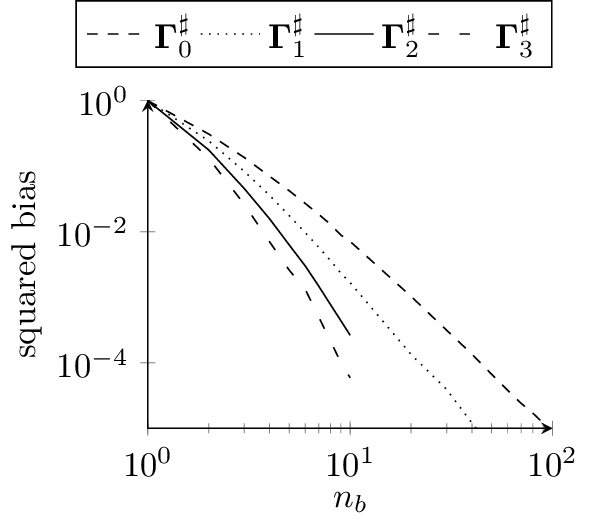}
&
\includegraphics[width=0.45\textwidth]{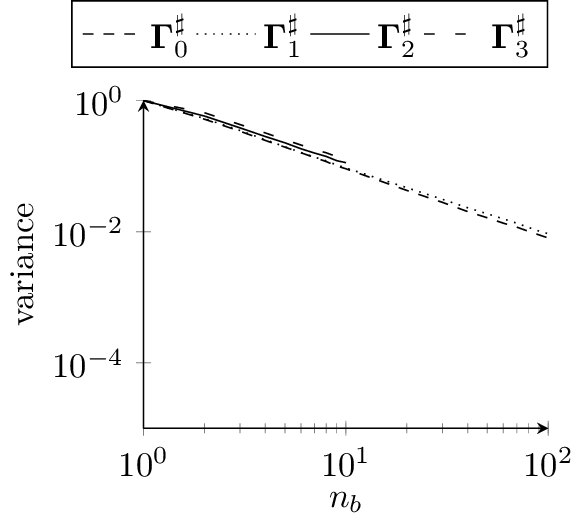}
\end{tabular}
\end{center}
\caption{\label{fig:SpeedTrees}Squared bias and variance for scenario trees}
\end{figure}

We are now able to compute the MSE using the protocol described in \S\ref{ssec:CompMSE}, for several values of the branching factor. We draw in Figure~\ref{fig:SpeedTrees} the evolution of both the squared bias and the variance with respect to the branching factor, for the strategy at each time step (here~$T=4$). For computational time reasons, we were not able to perform this experiment for strategies at the third and fourth time steps when the branching factor is more than~$10$ (the number of nodes involved at these time steps becomes too large).

On this figure which uses a logarithmic scaling on both axes, one can observe that the variance term on the right-hand side obviously dominates the squared bias term on the left-hand side. Moreover, while the convergence rate for the variance seems to be constant with respect to time, the bias decreases with a highest rate as time grows. This is natural since the tree has more nodes in the last time steps than in the first ones to represent the approximate strategy.

A careful inspection of the variance shows that the MSE appears experimentally to have a convergence rate of~$n_b^{-1}$ for every time step. In other words, the number of scenarios needed to obtain a given precision varies exponentially with respect to the time horizon: it is proportional to~$n_b^T$. Still, recall that in order to have a branching factor~$n_b$ equal to $10^2$, since we have a~$4$-period problem, we need to draw~$10^8$ scenarios!

\section{Particle methods} \label{sec:Particles}

We here propose a similar study using another numerical method based on scenarios: the particle method. This variational technique consists in first writing the optimality conditions for the problem. Then, scenarios are used to compute the expectations that lie in the optimality conditions. We briefly present the idea of this method and then study its error with respect to the number of scenarios used for sampling. A detailed description can be found in the work by~\citet{CarpentierCohenDallagi09}.

\subsection{Brief presentation}

The particle method consists in solving first-order optimality conditions that have been approximated through sampling. We first present the optimality conditions and then describe the sampling procedure. The resulting approximate optimality conditions are solved numerically using a gradient technique.

We suppose that the cost function and the system dynamics are differentiable with respect to both the state and the control, and that their derivatives are Lipschitz-continuous.

We suppose that Assumption~\ref{hyp:Indpdt} is still in force. In this context, \citet[Theorem 2.6]{CarpentierCohenDallagi09} state that, if a solution~$(\va{x}^*, \va{u}^*)$ to Problem~\ref{eqn:P} exists, then, for every time step $t=0, \dots, T$, there exists a random variable~$\va{\Lambda}_t^* \in L^2(\omeg, \trib, \prbt; \espacea{X})$, called the adjoint state, such that the following conditions hold:
\begin{subequations} \label{eqn:CondOptP}
	\begin{align}
		\va{x}_0^* &= \va{w}_0, \label{eqn:CondOptP-1} \\
		\va{x}_{t+1}^* &= f_t\left(\va{x}_t^*, \va{u}_t^*, \va{w}_{t+1}\right), \label{eqn:CondOptP-2} \\
		\va{\Lambda}_T^* &= \frac{\partial V}{\partial x}\left(\va{x}_T^*\right)^\top, \label{eqn:CondOptP-3} \\
		\va{\Lambda}_t^* &= \frac{\partial C_t}{\partial x}\left(\va{x}_t^*, \va{u}_t^*\right)^\top + \espcond{\frac{\partial f_t}{\partial x}\left(\va{x}_t^*, \va{u}_t^*, \va{w}_{t+1}\right)^\top \va{\Lambda}_{t+1}^*}{\va{x}_t^*}, \label{eqn:CondOptP-4} \\
		0 &= \frac{\partial C_t}{\partial u}\left(\va{x}_t^*, \va{u}_t^*\right)^\top + \espcond{\frac{\partial f_t}{\partial u}\left(\va{x}_t^*, \va{u}_t^*, \va{w}_{t+1}\right)^\top \va{\Lambda}_{t+1}^*}{\va{x}_t^*}. \label{eqn:CondOptP-5}
	\end{align}
\end{subequations}

Note that these optimality conditions are specific to the Markovian case and are not obtained in a straightforward manner. The authors first write first-order optimality conditions conditionally to the filtration~$\mathcal{F}_t$. Then Assumption~\ref{hyp:Indpdt} is used at each time step to recursively replace the conditional expectations with respect to~$\mathcal{F}_t$ by conditional expectations with respect to the optimal state~$\va{x}_t^*$.

Equations~\eqref{eqn:CondOptP-1} and~\eqref{eqn:CondOptP-2} simply recall the state dynamics of the system. Equations~\eqref{eqn:CondOptP-3} and~\eqref{eqn:CondOptP-4} introduce a backwards relation on the adjoint states. Actually, the adjoint state at time~$t$ may be seen as the sensitivity of the optimal cost from time~$t$ to the end of the time horizon with respect to the state variable at time~$t$. Hence there exists a relation between Equations~\eqref{eqn:CondOptP-3} and~\eqref{eqn:CondOptP-4} and the classical DP equation. Finally, Equation~\eqref{eqn:CondOptP-5} states that the derivative of the optimal future cost at a given time step with respect to the control at the same time step equals zero at the optimum. This is the classical first-order optimality condition when no additional constraints (apart for the non-anticipativity constraints) are imposed to the control.

Note that the adjoint state~$\va{\Lambda}_{t+1}^*$ in Equations~\eqref{eqn:CondOptP-4} and~\eqref{eqn:CondOptP-5} is, by construction, measurable with respect to the state variable~$\va{x}_{t+1}^* = f_t\left(\va{x}_t^*, \va{u}_t^*, \va{w}_{t+1}\right)$. Because of Assumption~\ref{hyp:Indpdt}, the conditional expectations in Equations \eqref{eqn:CondOptP-4} and \eqref{eqn:CondOptP-5} are actually expectations that are only supported by the random variable~$\va{w}_{t+1}$. One has to keep this in mind when discretizing these conditions.

We now want solve these conditions numerically using sampling techniques and gradient methods. Suppose we have $N$ samples $\va{w}^{\prime 1}, \dots, \va{w}^{\prime N}$, called scenarios or noise particles, for the random variable $\va{w}$. After the~$k$-th iteration, suppose we have $N$ samples for the current state $\va{x}^{(k)}$, control $\va{u}^{(k)}$, and adjoint state $\va{\Lambda}^{(k)}$. In order to compute conditions~\eqref{eqn:CondOptP-4} and~\eqref{eqn:CondOptP-5}, we need to evaluate expectations of functions involving~$\va{\Lambda}_{t+1}^{(k)}$, knowing the value~$\va{x}_t^{(k)}$ of the state at time $t$, for every particle~$k=1, \dots, N$. We have two ingredients to perform this operation:
\begin{enumerate}
	\item at the optimum, we know that for every time step $t=0, \dots, T-1$, the adjoint state $\va{\Lambda}_{t+1}^*$ is a function of the state $\va{x}_{t+1}^* = f_t(\va{x}_t^*, \va{u}_t^*, \va{w}_{t+1}$), where~$\va{u}_t^*$ is measurable with respect to~$\va{x}_t^*$, and~$\va{w}_{t+1}$ is independent of~$\va{x}_t^*$;
	\item we have $N$ samples of all the random variables involved, namely~$\va{x}^{(k)}$, $\va{u}^{(k)}$, $\va{\Lambda}^{(k)}$ and $\va{w}$.
\end{enumerate}
All we need is to define a regression operator, based on the current state and adjoint state samples, that associates with every value of the current state an estimate for the corresponding value of the adjoint state. At time $t$, we denote this regression operator by $\tilde{\va{\Lambda}}_t$. Using this material, we are now able to write the discretized optimality conditions, namely the initial condition on the state for every particle $i=1, \dots, N$:
\begin{subequations} \label{eqn:vitesse-CondOptPA}
	\begin{align}
		\va{x}_0^{\prime i} &= \va{w}_0^{\prime i}, \label{eqn:vitesse-CondOptPA-1} \\
	\intertext{the state dynamics for every time step $t=0, \dots, T-1$, and for every particle $i=1, \dots, N$:}
		\va{x}_{t+1}^{\prime i} &= f_t\left(\va{x}_t^{\prime i}, \va{u}_t^{\prime i}, \va{w}_{t+1}^{\prime i}\right), \label{eqn:vitesse-CondOptPA-2} \\
	\intertext{the final condition on the adjoint state, for every particle $i=1, \dots, N$:}
		\va{\Lambda}_T^{\prime i} &= \frac{\partial V}{\partial x}\left(\va{x}_T^{\prime i}\right)^\top, \label{eqn:vitesse-CondOptPA-3} \\
	\intertext{the backwards adjoint state dynamics, for every time step $t=T-1, \dots, 0$ and every particle $i=1, \dots, N$:}
		\va{\Lambda}_t^{\prime i} &= \frac{1}{N} \sum_{j=1}^N\Bigg( \frac{\partial C_t}{\partial x}\left(\va{x}_t^{\prime i}, \va{u}_t^{\prime i}\right)^\top \label{eqn:vitesse-CondOptPA-4} \\
		&\qquad + \frac{\partial f_t}{\partial x}\left(\va{x}_t^{\prime i}, \va{u}_t^{\prime i}, \va{w}_{t+1}^{\prime j}\right)^\top \tilde{\va{\Lambda}}_{t+1}\left(f_t\left(\va{x}_t^{\prime i}, \va{u}_t^{\prime i}, \va{w}_{t+1}^{\prime j}\right)\right) \Bigg), \nonumber \\
	\intertext{and the stationarity condition, for every time step $t=T-1, \dots, 0$ and every particle $i=1, \dots, N$:}
		0 &= \frac{1}{N} \sum_{j=1}^N \Bigg( \frac{\partial C_t}{\partial u}\left(\va{x}_t^{\prime i}, \va{u}_t^{\prime i}\right)^\top \label{eqn:vitesse-CondOptPA-5} \\
		&\qquad + \frac{\partial f_t}{\partial u}\left(\va{x}_t^{\prime i}, \va{u}_t^{\prime i}, \va{w}_{t+1}^{\prime j}\right)^\top \tilde{\va{\Lambda}}_{t+1}\left(f_t\left(\va{x}_t^{\prime i}, \va{u}_t^{\prime i}, \va{w}_{t+1}^{\prime j}\right)\right) \Bigg), \nonumber
	\end{align}
\end{subequations}
Note that in both equations~\eqref{eqn:vitesse-CondOptPA-4} and~\eqref{eqn:vitesse-CondOptPA-5}, the sums only involve the index~$j$, and hence the noise variable only. This corresponds to the fact that the conditional expectations in~\eqref{eqn:CondOptP-4} and~\eqref{eqn:CondOptP-5} with respect to~$\va{x}_t^*$ are in fact expectations.

The algorithm then consists in solving these conditions using a gradient method. Thus, each iteration $k$ of the particle method consists in three steps:
\begin{enumerate}
	\item integrate the $N$ state dynamics \eqref{eqn:vitesse-CondOptPA-1} and \eqref{eqn:vitesse-CondOptPA-2} using the current control particles $\va{u}^{\prime 1,(k)}, \dots, \va{u}^{\prime N,(k)}$;
	\item integrate the $N$ backwards adjoint state dynamics \eqref{eqn:vitesse-CondOptPA-3} and \eqref{eqn:vitesse-CondOptPA-4} using the current state and noise particles;
	\item compute the $N$ ``gradient particles'' $\va{g}^{\prime 1,(k)}, \dots, \va{g}^{\prime N,(k)}$ (the right-hand side of equation~\eqref{eqn:vitesse-CondOptPA-5}) and update every control particle~$i$ using a gradient step:
	\begin{equation*}
		\va{u}^{\prime i,(k+1)} = \va{u}^{\prime i,(k)} - \rho \va{g}^{\prime i, (k)}.
	\end{equation*}
\end{enumerate}
The algorithm stops when all gradient particles are sufficiently small, ensuring that condition \eqref{eqn:vitesse-CondOptPA-5} is almost fulfilled, for every~$i=1, \dots, N$.

\begin{remark}
	Note that, at each iteration, we have $N$ samples (or particles) for the state and control at every time step. Using a regression operator, we are able to compute a feedback function such as the one introduced in Section~\ref{ssec:MSE}, denoted by~$\va{\Gamma}^\sharp$. This feedback function associates a decision with every possible state of the system and this is the quantity on which we base our error analysis.
\end{remark}

\begin{figure}[htbp]
\begin{center}
\begin{tabular}{cc}
\includegraphics[width=0.45\textwidth]{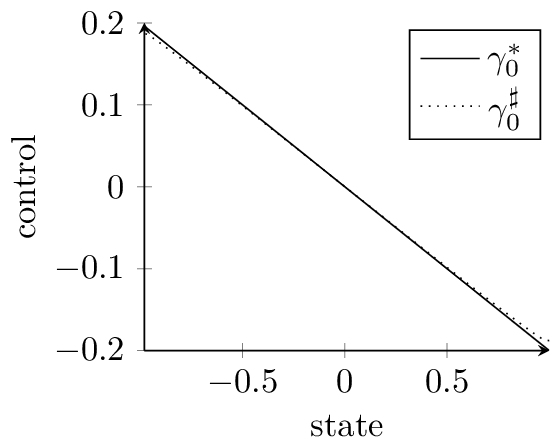}
&
\includegraphics[width=0.45\textwidth]{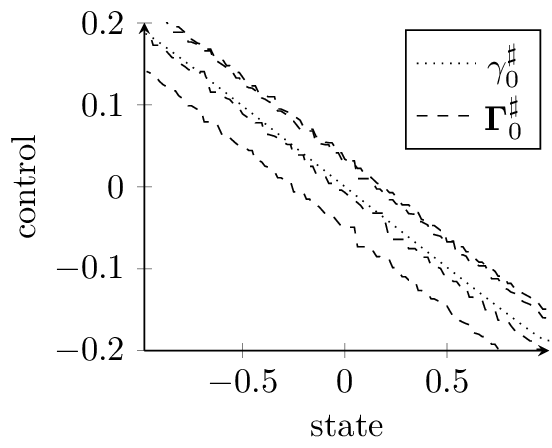}
\end{tabular}
\caption{\label{fig:StratExampleParticles}Strategies for particle methods}
\end{center}
\end{figure}

\subsection{Case study}

using the same example as in \S\ref{ssec:TreesCaseStudy}, we represent in the right-hand side of Figure~\ref{fig:StratExampleParticles} some samples of the approximate strategy (dashed curves) obtained by the particle method using~$3^4$ scenarios for discretization. Observe that these are piecewise constant functions with~$3^4$ pieces, instead of only $3$ pieces when using scenario trees. Indeed, using the particle method one has the same number of nodes at every time step. In other words, all the scenarios are used at each time step.
In the left-hand side of the same figure we draw the average approximate strategy (dotted curve) and the optimal strategy (solid curve). The bias term in the MSE is still the $L^2$-distance between the dotted curve and the solid curve, while the variance term is the average of the $L^2$-distances between the dashed curves and the dotted curve. Using the same number of scenarios as we did for for scenario trees in \S\ref{ssec:TreesCaseStudy}, we observe that both the bias and variance look much smaller to those observed in Figure~\ref{fig:StratExampleTrees}.

\begin{figure}[htbp]
\begin{center}
\begin{tabular}{cc}
\includegraphics[width=0.45\textwidth]{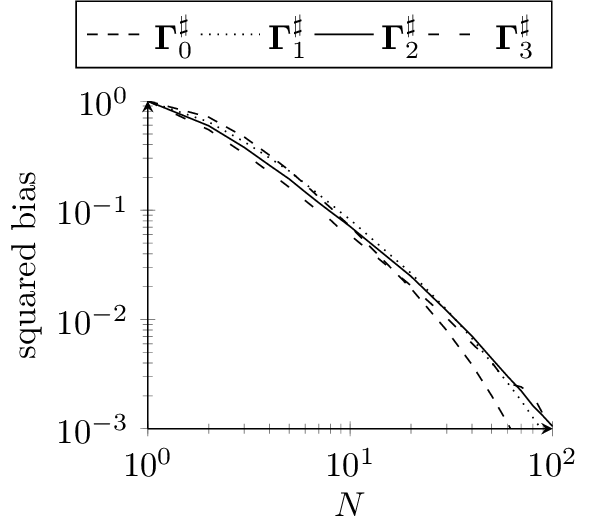}
&
\includegraphics[width=0.45\textwidth]{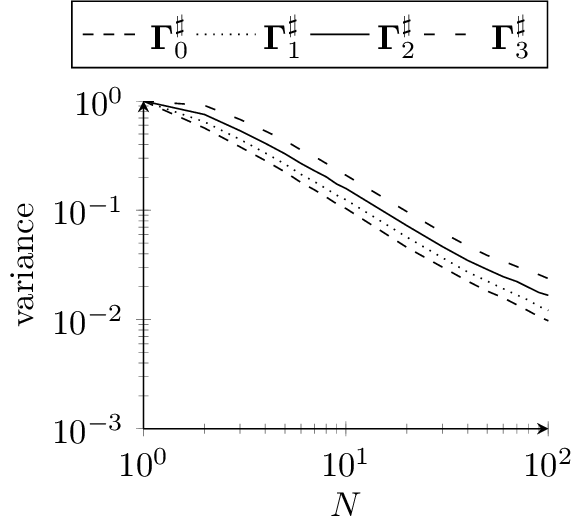}
\end{tabular}
\end{center}
\caption{\label{fig:SpeedParticles}Squared bias and variance for particles.}
\end{figure}

We can now compute the MSE for several values of the number~$N$ of scenarios used for discretization. We observe in figure~\ref{fig:SpeedParticles} how the squared bias and the variance associated with each strategy, i.e. for each time step, behave when the number of scenarios grows.
Note that the squared bias, on the left-hand side is dominated by the variance, on the right-hand side. A careful inspection of the latter shows that the term leads experimentally to a convergence rate that is slightly less than~$N^{-1}$. When comparing with the results obtained using scenario trees, one should notice that the $x$-axes do not refer to the same quantities in Figures~\ref{fig:SpeedTrees} and~\ref{fig:SpeedParticles}. Hence, we here observe a convergence rate in~$N^{-1}$, $N$ being the number of scenarios, while we observed a convergence rate of~$n_b^{-1} = N^{-\frac{1}{T}}$ in the case of scenario trees, which is much smaller.

Unlike scenario trees, the error associated with particle methods does not depend on the time horizon on this example: it is the same at every time step.

\section{The state space dimension issue} \label{sec:Dimension}

Most numerical methods that aim at solving stochastic optimal control problems encounter difficulties when the dimension of the state space becomes large. Since we are looking for strategies, that is functions that map a decision from every possible state of the system, even just storing such functions becomes challenging when the dimension of the state space is large. In the stochastic optimal control community framework, this difficulty is known as the \emph{curse of dimensionality}.

Note that particle methods are not prone to the same curse of dimensionality as DP, for instance. While DP encounters difficulties because it requires to build a somewhat regular grid to cover the state space, particle methods are designed to construct an adaptive mesh over the state space: the state particles are built so as to concentrate in regions that are often visited by the optimal state. As such, regardless of the dimension of the state space, whenever the dispersion of the optimal state is small, the particle method shall produce satisfying results.

\begin{figure}[htbp]
\begin{center}
\begin{tabular}{cc}
\includegraphics[width=0.45\textwidth]{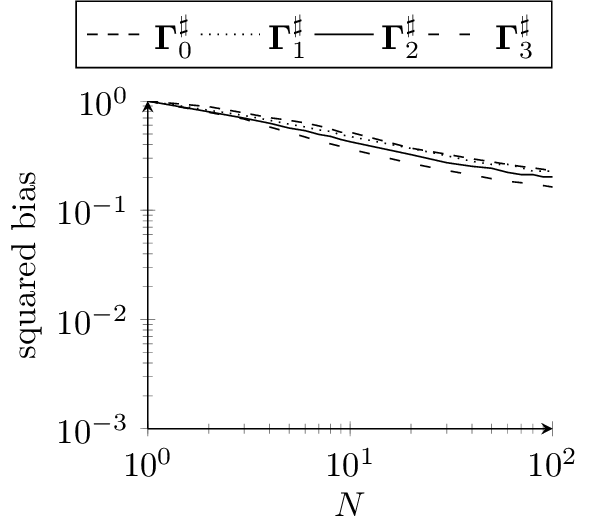}
&
\includegraphics[width=0.45\textwidth]{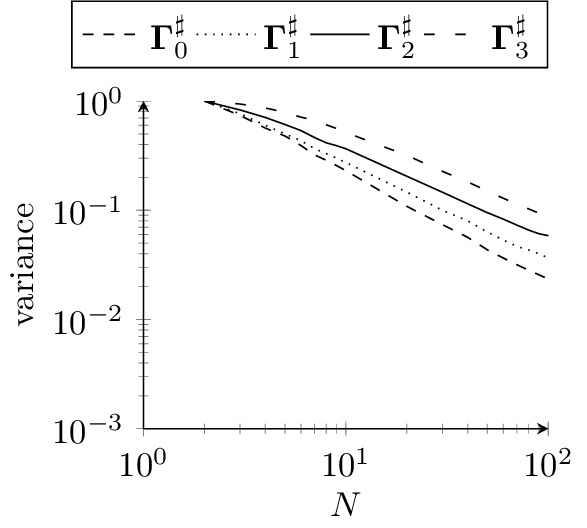}
\end{tabular}
\end{center}
\caption{\label{fig:SpeedParticles2d}Squared bias and variance for particles in dimension~2.}
\end{figure}

In the example we here study, the density of the optimal state does not seem to have this ``tightness property''. On a simple extension of Problem~\eqref{eqn:example} to a two-dimensional state space case, we observe in Figure~\ref{fig:SpeedParticles2d} how the squared bias and the variance associated with each strategy, i.e. for each time step, behave when the number of scenarios grows.
Note that the squared bias, on the left-hand side, now prevails over the variance term, on the right-hand side. A careful inspection of the latter shows that the term leads experimentally to a convergence that is slightly less than~$N^{-0.5}$. The variance itself also seems to decrease at a lower rate than in the one-dimensional case. Both of these observations indicate that the performance of the method could benefit from a better regression operator than the nearest neighbour that we used here.

Note that, for computational time reasons, we cannot produce such experimental results for higher dimensions. Hence we are not able to make a thorough study of the relation between the dimension of the state space and the convergence rate of both the squared bias and the variance.

On the other hand, when dealing with large state space dimensions, a clever idea is to make use of decomposition schemes in order to replace the solving of a high dimensional problem by the iterative solving of several smaller problems. This is indeed a common and powerful technique when using scenario trees~\citep[see][Ch. 3]{CarpentierCohenCulioli95,StochasticDecomposition96,StochasticProgramming03}. Both the improvement of the regression operator and the use of a decomposition scheme are the subject of ongoing research.

\section*{Conclusion}

We presented a numerical study concerning the comparison of two scenario-based approaches for solving stochastic optimal control problems with a discrete finite horizon. In the first section, we introduced the performance indicator that was used to compare both methods, namely the Mean Squared Error (MSE).

The first approach we considered was scenario tree modeling. We observed that the number of scenarios needed to obtain a given accuracy grew exponentially with the time horizon, making the implementation for multi-stage problems hardly tractable.

The second approach we considered was particle methods. In this case, we observed that the associated MSE does not depend on time. Thus, this recent approach seems to be well-suited for multi-stage problems. Unfortunately, it seems that its performance deteriorates rapidly when the dimension of the state space increases. However, this method is quite new and as such has received little attention. Its implementation and hence performance would certainly benefit from computational studies. Future works will be concerned with this topic.

\end{document}